\documentclass[lettersize,journal]{IEEEtran}
\usepackage{amsmath,amsfonts}
\usepackage{algorithm}
\usepackage{color}
\usepackage{algorithmic}
\usepackage{array}
\usepackage[caption=false,font=normalsize,labelfont=sf,textfont=sf]{subfig}
\usepackage{textcomp}
\usepackage{stfloats}
\usepackage{url}
\usepackage{verbatim}
\usepackage{graphicx}
\newtheorem{defn}{Definition}[section]
\newtheorem{rem}{Remark}[section]
\newtheorem{lem}{Lemma}[section]
\newtheorem{thm}{Theorem}[section]
\newtheorem{exm}{Example}[section]
\newtheorem{corollary}{Corollary}[section]
\newtheorem{prop}{Proposition}[section]

\def\BibTeX{{\rm B\kern-.05em{\sc i\kern-.025em b}\kern-.08em
    T\kern-.1667em\lower.7ex\hbox{E}\kern-.125emX}}
\begin{document}
\title{On reachability of Markov decision processes: a novel state-classification-based PI approach}
\author{Yanyun~Li, Xin~Guo$^{*}$
	and~Xianping~Guo
\thanks{Yanyun~Li is with School of Mathematics, Sun Yat-Sen University, Guangzhou, China (email: liyy536@mail2.sysu.edu.cn).}
\thanks{Xin~Guo is the corresponding author with the School of Science, Sun Yat-sen University, Guangzhou, China (email: guox87@mail.sysu.edu.cn).}
\thanks{Xianping~Guo is with the School of Mathematics Science, Sun Yat-Sen University; and Guangdong Province Key Laboratory of Computational Science,Guangzhou 510275, China (email: mcsgxp@mail.sysu.edu.cn).}}

\maketitle

\begin{abstract}
This paper concentrates on the reaching probability of discrete-time Markov decision processes (MDPs) with finite states and actions, and aims to give  a novel approach   for obtaining an optimal  policy that makes the MDPs have the minimal reaching probability. After establishing the existence of an optimal policy, in order to compute optimal  policies, we introduce the concept of an absorbing set for a stationary policy, and find   a computational method of the largest absorbing set. Using the largest absorbing set, we derive  an improved optimality equation (IOE), and prove  the uniqueness of solution of the IOE. By the uniqueness of solution of the IOE, we provide a novel  state-classification-based  policy iteration (PI) approach,  and prove that an optimal  policy and the minimal reaching probability can be obtained  by a new state-classification-based PI algorithm in a finite number of iterations. Moreover, we  prove that  the  number of calculations in one iteration of the PI algorithm here is less than that the PI algorithm for the  equivalent average  MDPs in  \cite{DM22}. Finally, an example in reliability and maintenance problems is given to illustrate our results.
\end{abstract}


\section{Introduction}
\IEEEPARstart{I}{n} this paper, we consider Markov decision processes (MDPs) in which, under a policy $\pi$, the state-action process    is denoted by $\{(Z_n, Y_n):\ n \geq 0\}$ on a $\pi$-dependent probability space  $(\Omega,\mathcal{F},P^{\pi})$. For any subset $B$ of a state space $S$ of the MDPs, let $\tau_B:=\inf\{n\geq 0:Z_n\in B\}$ be the life time of the MDPs. Then, for any initial state $i\in S$, the reaching probability of the MDPs  under a policy $\pi$ is defined as  $q^{\pi}_i:=P^{\pi}(\tau_B<\infty|Z_0=i)$. We want to minimize the reaching probability $q^{\pi}_i$ over the set $\Pi$ of all policies, and aim to give a novel efficient computational method for computing an optimal policy $\pi^*$ such that $q^{\pi^*}_i=\inf_{\pi\in\Pi}q^{\pi}_i$ for all $i\in S$.

Our study is inspired by MDPs
with risk neutral criteria	in \cite{BR11,OHLJBL:96,MLT05,S99} and the risk probability criterion in
\cite{BM04,HG-2020}, and is also motivated by the reachability problems in the probability of hitting a given set of states (e.g, \cite{DM22,DEJ11}). In MDPs with risk neutral criteria, the  existence and computational methods (e.g., the algorithms of policy iteration, linear programming and value iteration) of optimal policies have been deeply studied, see \cite{BR11,CR12,OHLJBL:96,MLT05,S99} for the risk-neutral discounted criterion and the risk neutral average criteria, and \cite{BR11,OHLJBL:96,MLT05,S99} for the risk neutral total reward criterion. As for the risk-sensitive case, the existence of an optimal policy is proven  in \cite{BR14,D07}, and a value iteration algorithm of an optimal policy is  given in \cite{BM04,HG-2020,Wen} for the risk probability criterion. As is well known, in the existing literature on MDPs with risk neutral/sensitive and probability criteria, the optimization criteria are defined by reward/cost functions. In reliability engineering situations \cite[p. 238-278]{MYK04}, there are different optimization objectives referring to reward/cost functions, for example, minimizing the deviation between the total system costs and the predetermined budgeted costs with uncertain random systems (e.g., \cite{XZ23}), and minimizing the long-run average cost of a non-repairable device whose degradation follows the Wiener process (e.g., \cite{YSX12}). However, there are some objectives not involving any reward/cost function, for instance, minimizing the reaching probability of series-parallel systems consisting of components which works independently with interchangeable redundancies, where the optimal policy is the order of the main components in ascending order of their reaching probabilities  (e.g., \cite{T23}).
Motivated from \cite{BR11,BM04,OHLJBL:96,HG-2020,MLT05,S99,T23}, this paper concentrates on the reachability problems in MDPs, instead of concerning about the components' orders.

As mentioned above, in this paper we mainly consider the reaching probability $q^{\pi}_i$ of MDPs on a set of failure  states $B$, which can be related not only to the failure probability in reliability engineering (e.g.,\cite{L96,N05,MYK04}), but also to the ruin probability in insurance-actuary \cite{AA10,TWWZ20}. We pay attention to
a computational method of an optimal policy that minimizes  the reaching probability  $q^{\pi}_i$. Different from the MDPs in \cite{BR11,CR12,OHLJBL:96,MLT05,S99,TSC92} with reward/cost functions, our optimization criterion  $q^{\pi}_i$ does not involve any reward/cost function. In particular, comparing with the risk-probability criterion in \cite{BM04,HG-2020,Wen}, which is based on reward functions and reward levels, our optimization criterion $q^{\pi}_i$  includes neither of them.
As far as we know, the reaching probability in discrete-time MDPs has been analyzed in \cite{DM22,DEJ11}, which is called the probability of reaching a  failure state set.  Chatterjee et al. \cite{DEJ11} discussed the minimal probability of reaching a failure state set while avoiding another set for discrete-time  MDPs with Borel state and action spaces, where it not only considered the existence of an optimal policy, but also connected its optimal control problem with a martingale characterization. For the reachability problem in \cite{DEJ11} with countable state and action spaces, $\acute{A}$vila et al. \cite{DM22} has transformed the problem of minimizing reaching probability into the equivalent (long-run) average MDPs by modifying the stochastic kernels and introducing a reward function, and then given the existence of an optimal policy and a linear programming algorithm of the minimal reaching probability and an optimal policy for finite states and actions by the results for MDPs in \cite{MLT05,S99}. As is well known, in addition to the linear programming algorithm, the minimal reaching probability and an optimal policy  in \cite{DM22} for finite states and actions can also be obtained by the policy iteration (PI) algorithm for the average MDPs in \cite{MLT05}. Each iteration of the PI algorithm for the  equivalent average  MDPs \cite{MLT05} consists of an improvement and evaluation step. The improvement step consists of two phases: First, improvement is sought through the first optimality equation. If no strict improvement is possible, a policy is established through the second optimality equation. Therefore, when the minimal reaching probability and an optimal policy  are computed by the PI algorithm for the equivalent average MDPs  in \cite{DM22}, each iteration in the PI algorithm from \cite{DM22}  needs to solve three equations.

Differing the PI algorithm for the equivalent average MDPs in  \cite{DM22}, in this paper we present a novel state-classification-based PI approach of computing the minimal reaching probability and optimal policies. Our new PI approach consists of  Algorithm 1 and Algorithm 2 below. Precisely, by  introducing  the concept of an absorbing set $F(g)$ for any stationary policy $g$ in Definition \ref{def4.1} below and analyzing the characterization of the absorbing set $F(g)$, we find a $B^c$-state-classification method (i.e., Algorithm 1), which can be used to obtain the  largest absorbing set $F_*$ in a finite number of iterations; see Proposition \ref{th4.4} below. For the special case of $F_*\subset \{\emptyset,B^c\}$,  we prove that an optimal policy has been obtained by Algorithm 1; see Remark \ref{re4.1} below. For the other case of  $F_*\not\subset \{\emptyset,B^c\}$, we derive  a so-called improved optimality equation (IOE), and prove  the uniqueness of a  solution of the IOE (i.e., Theorem \ref{le5.1}).  Using the uniqueness of the IOE, we provide a state-classification-based PI algorithm (i.e., Algorithm 2),
and proved that an optimal policy can be obtained by Algorithm 2 in a finite number of iterations (i.e., Proposition \ref{th5.3}). Comparison to the PI algorithm for the equivalent average MDPs in \cite{DM22} for solving the  minimal reaching probability,  the advantage of our Algorithm   2 is that each policy evaluation requires to solve only one equation (i.e., (\ref{5.2}) below) with $(|S|-|B|-|F_*|)$ variables, instead of solving three ones with $|S|$ variables in the PI algorithm for the equivalent average MDPs, where $|X|$ denotes the number of elements in any finite set $X$. Moreover, we prove that the  number of calculations (see its definition before Theorem \ref{re-5.1}) in one iteration of our PI approach (i.e., the sum of total number of calculations in Algorithm 1 and the number of calculations in one iteration of Algorithm 2), is more less than that in one iteration of the PI algorithm for the equivalent average MDPs in \cite{DM22}; see Theorem \ref{re-5.1} and Remark \ref{re-eff}.

The rest of our paper is organized as follows. We describe the MDP of MDPs in Section 2, and present some preparations in Section 3. In Section 4, we introduce the concept of an  absorbing set of a stationary policy and analyze its properties, and then give a $B^c$-state-classification method (i.e., Algorithm 1) of the largest absorbing set. In Section 5, by using the largest absorbing set, we set up the IOE with a unique solution, and present a state-classification-based PI
 algorithm (i.e., Algorithm 2) of optimal policies. Furthermore, Section 5 ends with presenting the comparison of the numbers of calculations of our state-classification-based PI approach (i.e., Algorithms 1-2) and the PI algorithm for the average MDPs in  \cite{DM22}; see Theorem \ref{re-5.1}  and Remark \ref{re-eff}. In Section 6, an example about a reliability and maintenance problem is given to illustrate  our state-classification-based PI approach.

\section{Problem Statement}\label{sec1}
The MDP considered in this paper is a four tuple:
\begin{equation}\label{2.1}
	\{S, B, (A(i): i \in S), p(\cdot| i, a)\},
\end{equation}
where the elements have the following explanations: $S$, a finite set, denotes the space of all states of the MDP under observation; $B \subset S$ is any given  subset  of $S$; $A(i)$  is a finite set of actions admissible at state $i \in S$; $p(\cdot|i, a)$, the transition probability on $S$ given $(i,a)\in \textbf{K}$, i.e., $p(j|i,a)\geq 0\ (j \in S)$ and $\sum_{j \in S} p(j | i,a)=1$, where  $\textbf{K}:=\{(i,a): i\in S, a\in A(i)\}$ denotes the set of all feasible state-actions pairs.

To describe the evolution of the process, let
\begin{equation}
	H_0:=S,\ \ H_n=\textbf{K}^{n}\times S,\ \ n\geq 1
\end{equation}
be the family of admissible histories up to time $n$.

\begin{defn}
	{\rm
		A randomized history-dependent policy $\pi=\{\pi_n:n\geq 0\}$ is a sequence of stochastic kernels~$\pi_n$~on~$A$~given~$H_n$ such that $\pi_n(A(i_n)|h_n)=1$b for all $n\geq 0$ and $h_n=(i_0, a_0,\cdots,i_{n-1},a_{n-1},i_n)\in H_n$, where $A:=\cup_{i \in S}A(i)$ is regarded as an  action space of the MDP. The set of all  randomized history-dependent policies is denoted by $\Pi$.
		
		A policy $\pi=\{\pi_n:n\geq 0\}\in \Pi$ is called a randomized Markov policy if, for each
		$h_n=(i_0,a_0,\cdots,i_{n-1},a_{n-1},i_n)$$\in H_n$ and $n\geq 0, \pi_n(\cdot |h_n)=:\pi_n(\cdot| i_n)$  depends only on  $i_n$.
		
		A randomized Markov policy $\pi=\{\pi_n:n\geq 0\} $ is called deterministic Markov if there exists a sequence of decision functions $\{g_n:n \geq 0\}$ such that for every $n \geq 0$, $g_n:i\mapsto g_n(i)\in A(i)$, and $\pi_n(g_n(i)|i)=1$ for all $i \in S$. We write such policy as $\pi=\{g_n:n \geq 0\}$.
		
		A deterministic Markov policy $\pi=\{g_n:n \geq 0\} $ is said to be stationary if $g_n=g$ for all $n \geq 0$. In this case, the policy $\pi$ is simply denoted by $g$. The set of all stationary  policies is denoted by $\Pi^s_d$.}
\end{defn}

Let $(\Omega, \mathcal{A})$~be the measurable space,~where~$\Omega:=H_{\infty}$~is the canonical sample space and~$\mathcal{A}$~is
the corresponding product~$\sigma$-algebra on $\Omega$. For any $i\in S$, $\pi=\{\pi_n:n\geq 0\} \in \Pi$,~by Tulcea's theorem in~\cite[p.147]{OHLJBL:96}, there exists a unique probability~$P^{\pi}_i$~on~$(\Omega, \mathcal{A})$ and a state-action process $\{(Z_n,Y_n):n\geq 0\}$ defined on $(\Omega, \mathcal{A})$ such that, for any $n\geq 0$, $j\in S$, $a \in A(i_n)$ and $h_n=(i_0, a_0,\cdots,i_{n-1},a_{n-1},i_n)\in H_n$, $P^{\pi}_i(Z_0 =i)=1$, $P^{\pi}_i(Y_n =a_n|h_n)=\pi_n(a_n| h_n)$ and $P^{\pi}_i(Z_{n+1}=j|h_n, Y_n=a_n)=p(j|i_n,a_n)$.

Now we present the evolution of the MDP (\ref{2.1}). Assume the initial state is $i_0\in S$. Given a policy $\pi=\{\pi_n: n\geq 0\}\in \Pi$, a controller selects an action $a_0\in A(i_0)$ with probability $\pi_0(a_0| i_0)$, and thus the process transfers to $i_1$ with probability $p(i_1|i_0, a_0)$ in one unit time. Then, during the next unit time, the controller chooses an action $a_1\in A(i_1)$ with probability $\pi_1(a_1|i_0, a_0,i_1)$, and the process transfers to $i_2$ with probability $p(i_2|i_1, a_1)$. The process goes on in this way.

To state the optimization problem here, let
\begin{eqnarray*}
	\tau_B:=\inf \{n \geq 0: \ Z_n \in B\}
\end{eqnarray*}
be the life time of the process (i.e., first hitting time on $B$).

For a given policy $\pi \in \Pi$ and initial state  $i\in S$ at time $n=0$, the reaching probability of the MDPs under policy $\pi$, denoted by $q_i^{\pi},$ is defined  by
\begin{eqnarray}\label{rp}
	q_i^{\pi}=P^{\pi}_i(\tau_B<\infty).
\end{eqnarray}
The reaching probability $q_i^{\pi}$ is called the failure probability in reliability engineering \cite{A07,L96,DH08,MYK04} (or the ruin probability in insurance-actuary \cite{AA10,TWWZ20}) (under policy $\pi$ and  initial state $i$).

Since $\tau_B=0$ for $i\in B$, it follows that $q^{\pi}_i=0\ (i\in B,\pi\in \Pi)$. Thus, for convenience, let
\begin{equation}\label{2.3}
	q_i^*:=\inf_{\pi \in \Pi}q_i^{\pi},\ \  i\in B^c,
\end{equation}
which  refers as the minimal reaching probability starting from state $i \in B^c$.
\begin{defn}\label{def2.3}
	\rm{
		A policy $\pi^*\in \Pi$ satisfying
		\begin{eqnarray*}
			q^{\pi^*}_i=q_i^*,\ \ i\in B^c,
		\end{eqnarray*}
		is said to be optimal.}
\end{defn}

The main goal of this paper is to introduce a new and efficient approach of finding the minimal reaching probability and its optimal policies.

\section{Preliminaries}\label{sec2}

In this section, we give the existence and property of an optimal policy.

The following lemma follows directly from the Markov property and Theorem 2.6 in \cite{Chen1992}.
\begin{lem}\label{le3.1}
\rm{	\begin{itemize}
		\item [(i)] For any $g \in \Pi_d^s$, $(q^g_i:i \in B^c)$ solves the following probabilistic equation system (PES):
		\begin{equation}\label{3.1}
			\begin{cases}
				x_i= p(B|i,g(i))+\sum\limits_{j\in B^c}p(j|i,g(i))x_j, & i\in B^c\\
				0\leq x_i\leq 1, \quad \quad \quad \quad \quad \quad \quad \quad \quad \quad \quad \quad & i\in B^c.
			\end{cases}
		\end{equation}
		\item [(ii)] If a nonnegative vector $(y_i:\ i\in B^c)$ satisfies
	\begin{eqnarray*}
		\begin{cases}
			y_i\geq p(B|i,g(i))+\sum\limits_{j\in B^c}p(j|i,g(i))y_j, &  i\in B^c\\
			y_i\geq 0, &  i\in B^c,
		\end{cases}
	\end{eqnarray*}
	then $y_i \geq q^g_i\ (i \in B^c)$.
	\end{itemize}}
\end{lem}

\begin{prop}\label{th3.1}	\begin{itemize}
\item[(i)] The minimal reaching probability $(q_i^*:\ i \in B^c)$ in $(\ref{2.3})$ satisfies the following optimality equation (OE):
	\begin{equation}\label{3.2}
		q^*_i=\min_{a\in A(i)}[p(B|i,a)+\sum_{j\in B^c} p(j|i,a) q^*_j],\quad i\in B^c.
	\end{equation}
	\item[(ii)]  There exists $g_* \in \Pi_d^s$ such that
	\begin{eqnarray*}
		&&p(B|i,g_*(i))+\sum_{j\in B^c} p(j|i,g_*(i))q^*_j\\
		&=&\min_{a\in A(i)}[p(B|i,a)+\sum_{j\in B^c}p(j|i,a)q_j^*],\ i \in B^c,
	\end{eqnarray*}
	and $g_*$ is optimal.
		\end{itemize}
\end{prop}
\begin{IEEEproof}
	(i) For any $\pi=\{\pi_n:\ n\geq 0\} \in \Pi$ and $i \in B^c$, by the definition of $q_i^{\pi}$,
	\begin{eqnarray*}
		q_i^{\pi}
		&=&\!\sum_{a \in A(i)}\!\!\pi_0(a|i)[p(B|i, a)+\sum_{j\in B^c}\!p(j|i,a)q_j^{\pi(i,a)}]\\
		&\geq&\min_{a\in A(i)}\{p(B|i, a)+\sum_{j\in B^c}(j|i,a)q_j^*\},
	\end{eqnarray*}
	where $\pi_n^{i,a}(\cdot|h_n):=\pi_{n+1}(\cdot|i,a,h_n)$ for all $h_n\in H_n$ and $n\geq 0$. Taking infimum over $\pi \in \Pi$ on the both sides,
	\begin{eqnarray*}
		q_i^* &\geq \min\limits_{a\in A(i)}[p(B|i, a)+\sum\limits_{j\in B^c}p(j|i,a)q_j^*], \ \ i\in B^c.
	\end{eqnarray*}
	Since $A(i)$ is finite, there exists $g_*\in \Pi^s_d$ such that
	\begin{eqnarray}\label{3.3}
		\nonumber	&&q_i^* \geq \min_{a\in A(i)}[p(B|i,a)+\sum_{j\in B^c}p(j|i,a)q_j^*]\\
		&=&p(B|i, g_*(i))+\sum_{j\in B^c}p(j|i,g_*(i))q_j^*, \ i\in B^c.
	\end{eqnarray}
	By Lemma~\ref{le3.1}(ii), $q^{g_*}_i\leq q^*_i\ (i\in B^c)$ and thus $q^{g_*}_i=q^*_i\ (i\in B^c)$. Hence, (i) follows from Lemma~\ref{le3.1}(i) and (\ref{3.3}). (ii) follows from (i) and Lemma \ref{le3.1}(ii).
\end{IEEEproof}

\begin{rem}\label{re3.2}  By Proposition \ref{th3.1}, we only need to pay our attention to the policies in $\Pi_d^s$ and to the minimal reaching probability $(q_i^*: i \in B^c)$ satisfying OE (\ref{3.2}). In the current case, $(q_i^*: i\in B^c)$ is the minimal nonnegative solution of the OE (\ref{3.2}) other than the unique one, see the following example for details.
\end{rem}
\begin{exm}\label{ex3.1}
	Consider a system with a state space $S=\{0,1,2\}$, $B=\{0\}$. Let $A(i):=\{c(i), d(i)\}\ (i \in S)$. The transition probabilities are given as follows:
	\begin{eqnarray*}
		p(j|i,c(i))=\delta_{ij},\ \ \
		p(j|i,d(i))=\delta_{0j}.
	\end{eqnarray*}	
	Consider the policies $g_1, g_2$ such that $g_1(i)=c(i)$ and $g_2(i)=d(i)$ for all $ i\in S$. Then, by Lemma \ref{le3.1}(i),
		$q_i^{g_2}=1$ for all $i=1,2$,
	i.e., $g_2$ satisfies OE (\ref{3.2}) but $q^{g_2}_1>q^{g_1}_1$ since $q^{g_1}_1=0$.
\end{exm}

Because of the non-uniqueness of the solution of OE  (\ref{3.2}), we next dedicate to IOE (\ref{3.2}) such that the IOE has a unique solution.

\section{On  $B^c$-absorbing sets of stationary policies}\label{sec3}

Since the solution of OE (\ref{3.2}) is not unique, for the computation of an optimal policy and the minimal reaching probability,  differing from that the PI algorithm for the equivalent average MDPs  in  \cite{DM22}, we will improve the OE (\ref{3.2}) to a system of equations with a unique solution.
For this purpose, we now introduce the concept of an absorbing set for  a stationary policy, which plays an important role in the following arguments.
\begin{defn}\label{def4.1}
	Let $g \in \Pi_d^s$. A set $F\subset B^c$ is called a $B^c$-closed set of $g$  if $p(F|i,g(i))=1\ (i \in F)$. The collection of all $B^c$-closed sets of $g$ is denoted by $\mathcal{C}_g$. $
	F(g):=\cup_{F\in \mathcal{C}_g}F
	$ is called the $B^c$-absorbing set of $g$.
\end{defn}

By Definition~\ref{def4.1}, $F(g)$ is actually the largest $B^c$-closed set of $g$. Denote $G(g):=B^c\setminus F(g)$.

The following lemma reveals that a better policy has a larger $B^c$-absorbing set.
\begin{lem}\label{le4.2}
\rm{	Suppose that $g,\tilde{g}\in\Pi_d^s$.
\begin{itemize}
	\item [(i)] $i\in F(g)$ if and only if $q^g_i=0$.
	
	\item [(ii)] If $q_i^{\tilde{g}}\leq q_i^g$ for all $i \in S$, then $F(g)\subset F(\tilde{g})$.
\end{itemize}	
}
\end{lem}

\begin{IEEEproof}
	We first prove (i). Obviously, $q^g_i=0\ (i\in F(g))$.
	Conversely, if $q^g_i=0$, then let $\tilde{F}=\{j \in S: \exists n\geq 0 \ {\rm s.t.} \ P^{g}_i(Z_n=j)>0\}$. Obviously, $\tilde{F}$ is closed. Note that $B$ can not be reached in finite steps from $i$. We know that $\tilde{F}$ is $B^c$-closed, and thus $i\in \tilde{F}\subset F(g)$.
	
	As for (ii), since $q_i^{\tilde{g}}\leq q_i^g$ for all $i \in S$, we have $q^{\tilde{g}}_i\leq q^g_i=0\ (i\in F(g))$. It follows from (i) that $F(g)\subset F(\tilde{g})$.
\end{IEEEproof}

Now we define
\begin{equation}\label{4.12}
	F_*:=\cup_{g \in \Pi_d^s}F(g), \quad G_*:=B^c\setminus F_*,
\end{equation}
where $F_*$ is called the largest absorbing set. Moreover, let
\begin{equation}\label{4.13}
	\Pi_d^s(F_*):=\{g\in \Pi_d^s:\ F(g)=F_*\}.
\end{equation}

The following lemma guarantees that $\Pi_d^s(F_*)\neq \emptyset$.
\begin{lem}\label{le4.3}
	Suppose that a policy $g_*\in\Pi_d^s$ is optimal. Then, $F(g_*)=F_*$.
\end{lem}

\begin{IEEEproof}
	It is obvious that $F(g_*) \subset F_*$ and thus we only need to prove that $F_* \subset F(g_*)$. Indeed, for any $g\in \Pi_d^s$, since $q_i^{g_*}\leq q_i^g$ for all $i \in S$, it follows from Lemma~\ref{le4.2}(ii) that $F(g)\subset F(g_*)$ and hence $F_* \subset  F(g_*)$.
\end{IEEEproof}
Next we consider the construction of $F_*$. For convenience, let $U_0:=B$ in the rest of this section, and
\begin{equation}\label{4.14}
	U_1:=\{i \in U_0^c:\ p(U_0|i,a)>0\quad \forall\ a \in A(i)\}.
\end{equation}
If $U_1=\emptyset$ or $U_0^c \setminus U_1=\emptyset$, then stop and define
\begin{equation}\label{4.15}
	A_{U_0}(i):=\{a \in A(i):\ p(U_0|i, a)=0\}
\end{equation}
for $i \in U_0^c \setminus U_1$ in the case $U_0^c\setminus U_1\neq \emptyset$. Recursively, for $n \geq 2$, if $U_{n-1}\neq\emptyset$ and $U_0^c\setminus \cup_{k=1}^{n-1}U_k\neq\emptyset$, then define
\small\begin{equation}\label{4.18}
	U_n\!:=\{i\in U_0^c\!\setminus\!\cup_{k=1}^{n-1}U_k:p( \cup_{k=0}^{n-1} U_k|i, a)>0\ \forall\ a \in A(i)\}.
\end{equation}
If $U_n=\emptyset$ or $U_0^c\setminus \cup_{k=1}^nU_k=\emptyset$, then stop and define
\begin{equation}\label{4.19}
	A_{U_{n-1}}(i):=\{a \in A(i):\ p(\cup_{k=0}^{n-1}U_k|i, a)=0\}
\end{equation}
for $i \in U_0^c \setminus \cup_{k=1}^{n}U_k$ in the case $U_0^c \setminus \cup_{k=1}^{n}U_k\neq \emptyset$.

\begin{thm}\label{th4.3}
\rm{	There exists $1 \leq N_* \leq |U_0^c|$ (the number of elements of $U_0^c$) such that
	\begin{itemize}
		\item [(i)]  $U_{N_*}= \emptyset$ or $U_0^c\setminus \cup_{k=1}^{N_*}U_k=\emptyset$;
		\item[(ii)] if $U_1=\emptyset$ or $U_0^c \setminus U_1=\emptyset$, then $N_*=1$; otherwise, $N_* \geq 2$, and for all $1 \leq n \leq N_*-1$,$
	U_{n}\neq \emptyset\ \ \text{and} \ \ U_0^c\setminus \cup_{k=1}^nU_k\neq \emptyset;$
	\item[(iii)] $F_*=U_0^c \setminus \cup_{k=1}^{N_*}U_k$. If $F_*\neq \emptyset$, then $A_{U_{N_*-1}}(i)\not= \emptyset\ (i\in F_*)$;
	 \item[(iv)]	if $F_*=\emptyset$, then $\Pi_d^s(F_*)=\Pi_d^s$; if $F_*\neq \emptyset$, then $\Pi_d^s(F_*)=\{g\in \Pi_d^s:\ g(i)\in A_{U_{N_*-1}}(i)\  \forall \  i\in F_*\}$.
     \end{itemize}}
\end{thm}
\begin{IEEEproof}
	It follows from (\ref{4.14}) and (\ref{4.18}) that
	$U_1, \cdots, U_{n}$ are disjoint. Let
	$
	N_*=\min\{N\geq 1: U_N=\emptyset\ {\rm{or}}\ U_0^c\setminus\cup_{k=1}^NU_k=\emptyset\}.
	$
	Obviously, $1\leq N_* \leq |U_0^c|$ since $U_0^c$ is finite. It is easy to see that (i) and (ii) hold for $N_*$.
	
	Next prove (iii). If $U_1=\emptyset$ or $U_0^c\setminus U_1=\emptyset$, then $N_*=1$. When $U_1=\emptyset$, we have $F_*\subset U_0^c=U_0^c\setminus U_1$. Conversely, since $U_0^c \setminus U_1=U_0^c \neq \emptyset$ and $U_0 \neq \emptyset$, we know that $A_{U_0}(i) \neq \emptyset$. Define $\tilde{g}\in \Pi_d^s$ such that
	\begin{equation}\label{4.22}
		\begin{cases}
			\tilde{g}(i)\in A(i),\quad & \text{if} \ i\in U_0,\\
			\tilde{g}(i)\in A_{U_0}(i),\quad &  \text{if}\ i\in U_0^c.
		\end{cases}	
	\end{equation}
	By (\ref{4.14})-(\ref{4.15}) and Definition \ref{def4.1}, we see that  $U_0^c\subset F(\tilde{g})\subset F_*$. Thus, $F_*=U_0^c=U_0^c\setminus U_1$. When $U_0^c\setminus U_1=\emptyset$, we have $F_*=\emptyset= U_0^c \setminus U_1$ from (\ref{4.12}),(\ref{4.14})-(\ref{4.15}) and Definition \ref{def4.1}.
	
	Now consider the case that $U_1\neq \emptyset$ and $U_0^c\setminus U_1\neq \emptyset$. In this case, $N_*\geq 2$. If $U_{N_*}\neq\emptyset$, then $U_0^c\setminus\cup_{k=1}^{N_*}U_k=\emptyset$. Therefore, we have $F(g)=\emptyset$ for every $g\in\Pi^s_d$. Indeed, suppose that $F(g)\neq \emptyset$. Take $i_0\in F(g)\subset U_0^c$. It can be proved that there exist $n(i_0)\leq N_*$, $k_1,\cdots,k_{n(i_0)}$ and $i_1\in U_{k_1}, i_2\in U_{k_2},\cdots, i_{n(i_0)}\in U_{k_{n(i_0)}}$ such that
	$k_1>k_2>\cdots >k_{n(i_0)-1}>k_{n(i_0)}=0$
	and
	\begin{eqnarray*}
		p(i_1|i_0,g(i_0))>0, \cdots, p(i_{n(i_0)}|i_{n(i_0)-1},g(i_{n(i_0)-1}))>0.
	\end{eqnarray*}
	Therefore,
	$q^{g}_{i_0}\geq p(i_1|i_0,g(i_0))\cdots p(i_{n(i_0)}|i_{n(i_0)-1},g(i_{n(i_0)-1}))>0$,
	which contradicts with $i_0\in F(g)$. Hence, $F_*=\emptyset=U_0^c\setminus\cup_{k=1}^{N_*}U_k$.
	
	While if $U_{N_*}=\emptyset$, then $U_0^c\setminus\cup_{k=1}^{N_*}U_k\neq \emptyset$. By Proposition~\ref{th3.1} and Lemma~\ref{le4.3}, there exists $g_*\in \Pi_d^s$ such that $F(g_*)=F_*$. Suppose that there exists $i_0\in F_*=F(g_*)$ but $i_0\notin U_0^c \setminus \cup_{k=1}^{N_*}U_k=U_0^c \setminus \cup_{k=1}^{N_*-1}U_k$. Then, it can also be proved that $q^{g_*}_{i_0}>0$,
	which contradicts with $i_0\in F(g_*)$.
	Therefore, $F_*=F(g_*)\subset U_0^c \setminus \cup_{k=1}^{N_*}U_k$. On the other hand, since $U_0^c \setminus \cup_{k=1}^{N_*}U_k\neq \emptyset$ and $\cup_{k=0}^{N_*-1}U_k \neq \emptyset$, we obtain $A_{U_{N_*-1}}(i) \neq \emptyset$. Take $\tilde{g}\in \Pi_d^s$ such that
	\begin{equation}\label{4.23}
		\begin{cases}
			\tilde{g}(i)\in A(i),\quad & if\ i\in \cup_{k=0}^{N_*}U_k\\
			\tilde{g}(i)\in A_{U_{N_*-1}}(i),\quad & if\ i\in U_0^c \setminus \cup_{k=1}^{N_*}U_k.
		\end{cases}
	\end{equation}
	
	From (\ref{4.12}), (\ref{4.15}) and (\ref{4.19}), we obtain $U_0^c \setminus \cup_{k=1}^{N_*}U_k\subset F(\tilde{g})\subset F_*$. (iii) is proved.
	
	Finally, (iv) follows from Lemma \ref{le4.3} and Proposition \ref{th3.1}.
\end{IEEEproof}
By Lemma~\ref{le4.3}, we only need to consider the policies in $\Pi_d^s(F_*)$. Furthermore, from Theorem~\ref{th4.3} and its proof, we see that $\Pi_d^s(F_*)$ can be obtained by collecting policies satisfying (\ref{4.23}) in the case $N_* \geq 2$ or (\ref{4.22}) in the case $N_*=1$.

For any $g\in\Pi_d^s(F_*)$, define
\begin{equation}\label{4.1}
	G_1(g):=\{i \in U_0^c:\ p(U_0|i, g(i))>0\}.
\end{equation}
Recurrently, for $n\geq 1$, if $U_0^c \setminus \cup_{k=1}^{n}G_k(g)\neq \emptyset$, then define
\begin{equation}\label{4.2}
	G_{n+1}(g):=\{i \in U_0^c \setminus \cup_{k=1}^{n}G_k(g):\ p(G_{n}(g)|i, g(i))>0\}.
\end{equation}
Note that in the definition of $G_{n+1}(g)$, $\ p(G_{n}(g)|i, g(i))>0$ is equivalent to $\ p(\cup_{k=0}^nG_{k}(g)|i, g(i))>0$ (where $G_0(g)=U_0$). By Theorem~\ref{th4.3}, we have the following result.
\begin{corollary}\label{le4.1}
\rm{	For any $g \in \Pi_d^s(F_*)$, there exists a unique $1\leq N_g \leq |U_0^c|$ such that
	\begin{itemize}
		\item [(i)] $G_{N_g}(g)=\emptyset$ or $U_0^c\setminus \cup_{k=1}^{N_g}G_k(g)=\emptyset$;
		
		\item[(ii)] if $G_1(g)=\emptyset$ or $U_0^c \setminus G_1(1)=\emptyset$, then $N_g=1$; otherwise, $N_g \geq 2$, and for all $1 \leq n \leq N_g-1$, we have
		$G_{n}(g)\neq \emptyset$ and $ U_0^c\setminus \cup_{k=1}^nG_k(g)\neq \emptyset$;
		\item[(iii)]  $F_*=F(g)=U_0^c \setminus \cup_{k=1}^{N_g}G_k(g)$.
	\end{itemize}}
\end{corollary}
The following Theorem~\ref{th4.2} reveals that for every $g\in \Pi_d^s(F_*)$, $(q_i^g:\ i \in  G_*)$ (if $G_*\neq \emptyset$)  is the unique solution of
\begin{equation}\label{4.7}
	\begin{cases}
		x_i= p(B|i,g(i))+\sum\limits_{j\in G_*}p(j|i,g(i))x_j, &  i\in G_*\\
		0\leq x_i\leq 1, \quad \quad \quad \quad \quad \quad \quad \quad \quad \quad \quad \quad  & i\in G_*,
	\end{cases}
\end{equation}	
which improves (\ref{3.1}) and will play a key role in the PI algorithm.

\begin{thm}\label{th4.2}
\rm{	Suppose that $g \in \Pi_d^s(F_*)$.
	\begin{itemize}
		\item [(i)] If $F_*=U_0^c$, then $q^g_i=0$ for all $i\in U_0^c$.
		
		\item[(ii)] If $F_*=\emptyset$, then $q_i^g=1$ for all $i\in U_0^c$.
		
		\item[(iii)]  If $G_*\not=\emptyset$, then  $(q_i^g:\ i \in G_*)$ is the unique solution of $(\ref{4.7})$.
	\end{itemize}
}
\end{thm}
\begin{IEEEproof}
	(i) follows from (\ref{4.12}) and Lemma \ref{le4.2}. We now prove (ii). If $F_*=\emptyset$, then $U_0^c \setminus \cup_{k=1}^{N_g}G_k(g)= \emptyset$. By Corollary \ref{le4.1}(i),
	\begin{eqnarray*}
		G_k(g)\neq \emptyset,\ \ \ k=1,\cdots, N_g.
	\end{eqnarray*}
	Denote
	$
	\delta=\min_{i\in G_k(g), 1\leq k\leq N_g}\limits p(G_{k-1}(g)|i,g(i)),
	$
	where $G_0(g)=U_0$. By the finiteness of $G_k(g)\ (1\leq k\leq N_g)$, we know that $\delta>0$.
	
	It can be proved that for all $k\in \{1,\cdots,N_g\},\ i\in G_k(g)$,
	\begin{equation}\label{4.8}
		P_i^g(\tau_{_{U_0}}\leq k+1)\geq P_i^g(\tau_{_{U_0}}= k+1)\geq \delta^{k+1}.
	\end{equation}
	Indeed, $P_i^g(\tau_{_{U_0}}\leq 1)\geq P_i^g(\tau_{_{U_0}}=1)=p(B|i,g(i))\geq \delta$ for any $i\in G_1(g)$.
	Recursively, for any $i\in G_k(g)\ (2\leq k\leq N_g)$,
	\begin{eqnarray*}
		&&P_i^g(\tau_{_{U_0}}=k+1)=\sum_{j\in U_0^c}p(j|i,g(i))P_j(\tau_{_{U_0}}=k)\\
		&\geq& \delta^k p(G_{k-1}(g)|i,g(i))\geq \delta^{k+1}.
	\end{eqnarray*}
	
	Now, since $U_0^c \setminus \cup_{k=1}^{N_g}G_k(g)= \emptyset$, we know that $P_i^g(\tau_{_{U_0}}> N_g+1)\leq 1-\delta^{N_g+1}$ for any $i\in U_0^c=\cup_{k=1}^{N_g}G_k(g)$.
	
	Furthermore, by the Markov property and mathematical induction, it can be proved that for any $i\in U_0^c$,
	\begin{eqnarray*}
		P_i^g(\tau_{_{U_0}}> k(N_g+1))\leq (1-\delta^{N_g+1})^k \ \ \forall\ k\geq 1.
	\end{eqnarray*}
	Letting $k\uparrow\infty$ yields $q^g_i=P_i^g(\tau_{_{U_0}}<\infty)=1$. (ii) is proved.
	
	Now turn to prove (iii). By Theorem 2.2 in \cite{Chen1992}, (\ref{4.7}) has a solution, denoted as $(\tilde{q}_i:i\in G_*)$.  Define $\tilde{y}_i:=\tilde{q}_i$ for $i\in G_*$ and $0$ for $i\in F_*$.
	Obviously, $(\tilde{y}_i:i\in U_0^c)$ is a nonnegative solution of the PES (\ref{3.1}). By Lemma \ref{le3.1}, $(q^g_i:i\in U_0^c)$ is the nonnegative minimal solution, we have $q^g_i=0\ (i\in F_*)$ and hence $(q_i^g:i \in G_*)$ is a solution of $(\ref{4.7})$.
	
	Next prove that $(q_i^g: i \in G_*)$ is the unique solution of $(\ref{4.7})$. To do so, we rewrite (\ref{4.7}) as $
		[I-(\emph{\textbf{P}}^g_{G_*})^n] \emph{\textbf{x}}=V^g_{U_0}+ \cdots+(\emph{\textbf{P}}^g_{G_*})^{n-1}V^g_{U_0}$,
	where $ V^g_{U_0}=(p(U_0|i,g(i)): i \in U_0^c)^T$,
	$\emph{\textbf{P}}^g_{G_*}=(p(j|i,g(i)):i,j \in G_*)$, $I$ is the identity matrix and $\emph{\textbf{x}}=(x_k: k \in G_* )^T$.
	
	It suffices to prove that there exists $n<\infty$ such that $[I-(\emph{\textbf{P}}^g_{G_*})^n]$ is invertible.
	For this aim, we first verify that there exists $n <\infty$ such that for all $i \in G_*$,
	\begin{equation}\label{4.11}
		\sum_{j \in G_*}p^{(n)}(j|i, g(i))<1,
	\end{equation}
	where $p^{(n)}(j|i, g(i))$ is the $(i,j)$'th element of $(\emph{\textbf{P}}^g_{G_*})^n$. In fact, if $i \in G_1(g)$, then by (\ref{4.1}),
	$\sum_{j \in G_*}p(j|i, g(i))\leq \sum_{j \in U_0^c}p(j|i, g(i))<1.
	$
	If $i \in G_k(g)$ for some $k\in \{2,\cdots,N_g\}$, then by (\ref{4.2}), there exist $i_{k-1}\in G_{k-1}(g), \cdots, i_2\in G_2(g), i_1\in G_1(g)$ such that
	$
	p(i_{k-1}|i,g(i))>0, \cdots, p(i_{1}|i_2,g(i_2))>0.
	$
	By (\ref{4.1}), we have $\sum_{j\in U_0}p(j|i_1,g(i_1))>0$.
	
	 Hence, $\sum_{j \in U_0} p(j|i_1, g(i_1)) \cdots p(i_{k-2}|i_{k-1}, g(i))p(i_{k-1}|i, g(i))>0$,
	which implies that $\sum_{j \in U_0}p^{(k)}(j|i, g(i))>0$, and thus $\sum_{j \in G_*}p^{(k)}(j|i, g(i))\leq \sum_{j \in U_0^c}p^{(k)}(j|i,g(i))<1$.
	Therefore, for any $i \in G_*$, there exists $n(i)$ such that $\sum_{j \in G_*}p^{(n(i))}(j|i, g(i))<1$.
	Furthermore, for any $i\in G_*$ and $n\geq n(i)$, $\sum_{j \in G_*}p^{(n)}(j|i, g(i))\leq \sum_{l\in G_*}p^{(n(i))}(l|i, g(i))<1$. Since $G_*$ is finite, we know that (\ref{4.11}) holds for $n=\max_{i\in G_*}n(i)$. Hence, $I-(\emph{\textbf{P}}^g_{G_*})^n$ is invertible.
\end{IEEEproof}

By Theorem~\ref{th4.3}, we have an algorithm of $F_*$ and $\Pi_d^s(F_*)$ below.

\quad

{\bf Algorithm 1:}{ A $B^c$-state-classification method.}\label{al4.2}

\hspace*{0.02in}{\bf Input:}
The failure state set $B$ and the transition probabilities $p(j|i, g(i))$ for every $i,j \in B^c$ and $g \in \Pi^s_d$.

\hspace*{0.02in}{\bf Output:}
The largest $B^c$-absorbing set $F_*$ and the new policy set $\Pi_d^s(F_*)$.

\hspace*{0.02in}{\bf Step 1:} Set $n=1$ and $U_0:=B$, obtain $U_1$ by (\ref{4.14}).

 A.
 If $U_1=\emptyset$, then stop, set $F_*=U_0^c$, $N_*=1$ and get $A_{U_0}(i)$ by (\ref{4.19}) for $i \in F_*$. Furthermore, set
    \begin{eqnarray*}
    	\Pi_d^s(F_*)=\{g\in \Pi_d^s:\ g(i)\in A(i)\ \text{for}\ i\in U_0\\
    	\text{and}
    	\ g(i)\in A_{U_0}(i)\ \text{for} \ i\in F_*\};
\end{eqnarray*}

 B.
 if $U_0^c\setminus U_1=\emptyset$, then stop and set $F_*=\emptyset$, $N_*=1$ and $\Pi_d^s(F_*)=\Pi_d^s$;

 C.
 otherwise increment $n$ by $1$ and go to step $2$.

 \hspace*{0.02in}{\bf Step 2:} Collect states $i \in U_0^c \setminus \cup_{k=1}^{n-1}U_k$ satisfying $p(\cup_{k=0}^{n-1}U_k|i, a)>0$ for all $a \in A(i)$, and then obtain $U_n$.

 \hspace*{0.02in}{\bf Step 3:} A.
 If $U_n=\emptyset$, then stop, set $F_*=U_0^c \setminus \cup_{k=1}^{n-1}U_k$, $N_*=n$ and get $A_{U_{n-1}}(i)$ by (\ref{4.19}) for $i \in F_*$. Then, set
 \begin{eqnarray*}
 	\Pi_d^s(F_*)=\{g\in \Pi_d^s:\ g(i)\in A(i)\ \text{for}\ i\in S\setminus F_*\\ \ \text{and}\ g(i)\in A_{U_{n-1}}(i)\ \text{for}\ i\in F_*\};
 \end{eqnarray*}

 B.
 if $U_0^c\setminus \cup_{k=1}^nU_k=\emptyset$, then stop, set $F_*=\emptyset$, $N_*=n$ and $\Pi_d^s(F_*)=\Pi_d^s$;

  C.
 otherwise increment $n$ by $1$ and turn back to step $2$.

\quad

\begin{prop}\label{th4.4} The  $F_*$ and $\Pi_d^s(F_*)$  can be obtained by Algorithm 1 in a finite number of iterations.
\end{prop}
\begin{IEEEproof}
	Obviously, it follows from Theorem~\ref{th4.3} and the description of Algorithm 1.
\end{IEEEproof}
\begin{rem}\label{re4.1}
	By Algorithm 1, if $F_*=\emptyset$, then $\Pi_d^s(F_*)=\Pi_d^s$ and by Theorem~\ref{th4.2}(ii), we know that $q^g_i=1\ (i\in B^c)$ for all $g\in \Pi_d^s$, and hence $q^*_i=1$ for all $i\in B^c$ and every policy $g\in \Pi_d^s$ is optimal. While if $G_*=\emptyset$ (i.e., $F_*=B^c$), then it follows from Theorem~\ref{th4.2}(i) that  $q^g_i=0\ (i\in B^c)$ for all $g\in \Pi_d^s(F_*)$, and hence $q^*_i=0$ for all $i\in B^c$ and each policy $g\in \Pi_d^s(F_*)$ is optimal.
\end{rem}

\section{A state-classification-based PI algorithm of optimal policies}\label{sec4}

By Remark~\ref{re4.1}, if $F_*=\emptyset$ or $G_*=\emptyset$, then the optimal policy and the minimal reaching probability have been obtained in the previous section. Therefore, it remains to consider the case that $F_*\neq \emptyset$ and $F_*\neq B^c$ (and hence $G_*\not=\emptyset$) in this section. Since an optimal policy must be in $\Pi_d^s(F_*)$, we mainly consider how to find an optimal policy  in $\Pi_d^s(F_*)$ and obtain the minimal reaching probability $(q^*_i : i\in B^c)$.

\begin{thm}\label{le5.1}A policy $g$ in $\Pi_d^s(F_*)$ is optimal if and only if $(q^g_i: i\in G_*)$ solves the following IOE (\ref{im-OE1}):
	\begin{eqnarray}\label{im-OE1}
		u_i=\min_{a\in A(i)}\{p(B|i,a)+\sum_{j\in G_*} p(j|i,a) u_j\},\ \ i\in G_*.
	\end{eqnarray}
\end{thm}
\begin{IEEEproof}
	Suppose that $g\in\Pi_d^s(F_*)$ is optimal. Then, $q^*_i=q^g_i\ (i\in B^c)$ and $q^g_j=0\ (j\in F_*)$ (by Lemmas \ref{le4.2}-\ref{le4.3}). Thus, by Proposition \ref{th3.1}  we have
	\begin{eqnarray}\label{5.1}
		\nonumber	q^g_i\!\!&=&\!q^*_i=\min_{a\in A(i)}\{ p(B|i,a)+\sum_{j\in B^c} p(j|i,a) q^*_j\}\\
		&=&\min_{a\in A(i)}\{ p(B|i,a)+\sum_{j\in G_*} p(j|i,a) q^{g}_j\},\ \ i\in G_*,
	\end{eqnarray}
	i.e., $(q^g_i: i\in G_*)$ solves the IOE (\ref{im-OE1}).
	
	Given $g \in \Pi_d^s(F_*)$, suppose that $(q^g_i:i\in G_*)$ solves (\ref{im-OE1}). Then, by (\ref{4.13}) we have $F(g)=F_*$ and hence $G(g)=G_*$. For the policy $g_*$ in Proposition \ref{th3.1}(ii), by Theorem~\ref{th4.2} and Lemma \ref{le4.3}, we see that $q^g_i=q^{g_*}_i=0$ for all $i\in F_*$, and
	\begin{eqnarray*}
		q^g_i&=&\min_{a\in A(i)}\{ p(B|i,a)+\sum_{j\in G_*} p(j|i,a) q^{g}_j\}\\
		&\leq& p(B|i,g_*(i))+\sum_{j\in G_*} p(j|i,g_*(i)) q^{g}_j,\ i\in G_*.
	\end{eqnarray*}
	On the other hand, since  $g_*$ is optimal and $q^{g_*}_i=0\ (i\in F_*)$, by (\ref{3.2}) we have
	\begin{eqnarray*}
		q^{g_*}_i&=&\min_{a\in A(i)}\{ p(B|i,a)+\sum_{j\in G_*} p(j|i,a) q^{g_*}_j\}\\
		&=&p(B|i,g_*(i))+\sum_{j\in G_*} p(j|i,g_*(i)) q^{g_*}_j, i\in G_*.
	\end{eqnarray*}
	Therefore, for any $i\in G_*$,
	$0\leq q^g_i-q^{g_*}_i\leq \sum\limits_{j\in G_*} p(j|i,g_*(i))(q^g_j- q^{g_*}_j).$
	Recursively,
	\begin{eqnarray*}
		0\leq q^g_i-q^{g_*}_i\leq \sum_{j\in G_*} p^{(n)}(j|i,g_*(i))(q^g_j- q^{g_*}_j),
	\end{eqnarray*}
	where $p^{(n)}(j|i, g_*(i))$ is the $(i,j)$'th element of $(\emph{\textbf{P}}^{g_*}_{G_*})^n$ and $\emph{\textbf{P}}^{g_*}_{G_*}=(p(j|i,g_*(i)): i,j\in G_*)$.
	Since $G_*$ is finite, it follows from the above inequality and a same argument as in the proof of (\ref{4.11}) that $q^g_i=q^{g_*}_i=q^*_i$ for all $i\in B^c=F_*\cup G_*$.
	Hence, $g$ is optimal.
\end{IEEEproof}

To compare two policies in $\Pi_d^s(F_*)$, we define $A^*(i)$ for $i\in S$:
\begin{equation}\label{4.20}
	A^*(i):=\begin{cases}
		A(i),& \text{if} \ i\in B\cup G_*\\
		A_{U_{N_*-1}}(i), & \text{if} \ i\in F_*,
	\end{cases}
\end{equation}
which are not empty (by Theorem \ref{th4.3}(i)).

\begin{thm}\label{th5.1}
\rm{	Given $g\in \Pi_d^s(F_*)$. For $i \in B^c$, let
	\begin{eqnarray*}
		A_g(i):=\{a \in A^*(i):\ q_i^g>p(B|i,a)+\sum_{j \in G_*}p(j|i, a)q_j^g\}.
	\end{eqnarray*}
	Define a policy $\tilde g$  as follows: for $i\in B^c$,
\begin{itemize}
	\item [(i)] if $A_g(i)=\emptyset$, let  $\tilde{g}(i):=g(i)$;
	
	\item [(ii)] if $A_g(i)\neq \emptyset$, take 	
	\begin{eqnarray*}
		\tilde{g}(i)\in\mathop{\rm{argmin}}_{a\in A_g(i)}\{p(B|i,a)+\sum_{j \in G_*}p(j|i, a)q_j^g\}.	
	\end{eqnarray*}
\end{itemize}	
	Then $q_i^{\tilde{g}}\leq q_i^g (i\in B^c)$. Moreover, if $\tilde{g} \neq g$, then $q_i^{\tilde{g}}<q_i^g$ for some $i \in B^c$.}
\end{thm}
\begin{IEEEproof}Since $q^g_i=0 \ (i\in F_*)$, for any $i \in B^c$, if $A_g(i) \neq \emptyset$, we have
	$q_i^g>p(B|i,\tilde{g}(i))+\sum_{j \in B^c}p(j|i, \tilde{g}(i))q_j^g.$ If $A_g(i)=\emptyset$ then $g(i)=\tilde{g}(i)$ and $
	q_i^{g}=p(B|i,\tilde{g}(i))+\sum_{j\in B^c}p(j|i,\tilde{g}(i))q_j^g.
	$
	Thus, for $i \in B^c$,
	$q_i^{g}\geq p(B|i,\tilde{g}(i))+\sum_{j \in B^c}p(j|i,\tilde{g}(i))q_j^{g}.$
	Therefore, by Lemma \ref{le3.1}, we have $q_i^g \geq q_i^{\tilde{g}}$ for all $i \in B^c$.
	
	If $\tilde{g}\neq g$, then $A_g(i) \neq \emptyset$ and $\tilde{g}(i)\in A_g(i)$ for some $i\in B^c$. By Theorem~\ref{th4.2}, we can prove that $q_i^g>q_i^{\tilde{g}}$ for some $i\in B^c$.
\end{IEEEproof}

We now present a state-classification-based PI algorithm for computing optimal policies in $\Pi_d^s(F_*)$. By Theorem~\ref{le5.1} and Theorem~\ref{th4.2}(iii), for any $g\in \Pi_d^s(F_*)$, we have $q^g_i=0\ (i\in F_*)$  and $(q^g_i: i\in G_*)$  can be obtained by solving (\ref{5.2})
\begin{eqnarray}\label{5.2}
	x_i=p(B|i, g(i))+\sum_{j \in G_*}p(j|i, g(i))x_j,\ \ i \in G_*,
\end{eqnarray}
which, together with Theorem \ref{th4.2}(iii), implies that
\begin{eqnarray*}
	Q^g_{G_*}=(I-\emph{\textbf{P}}^g_{G_*})^{-1}V_B^g,
\end{eqnarray*}
where $Q^g_{G_*}=(q_i^g:i \in G_*)^T, \  \emph{\textbf{P}}^g_{G_*}=(p(j|i, g(i)) : \ i,j \in G_*)$ and $ V_B^g=(p(B|i, g(i)):i \in G_*)^T$.

\quad

{\bf Algorithm 2:} A state-classification-based PI algorithm.

\hspace*{0.02in}{\bf Input:} The failure state set $B$ and the transition probabilities $p(j|i, g(i))$ for every $i,j \in B^c$ and $g \in \Pi^s_d$.
	
	\hspace*{0.02in}{\bf Output:} The optimal policy $g_*$ and the minimal reaching probability $(q_i^*:\ i \in B^c)$.
	
	\hspace*{0.02in}{\bf Step 1:} Obtain $F_*$ and $\Pi_d^s(F_*)$ by Algorithm 1.
		
	\hspace*{0.02in}{\bf Step 2:} Set $n=0$, select an arbitrary policy $g_0 \in \Pi_d^s(F_*)$, i.e., $g_0(i) \in A^*(i),\ i\in B^c$.
		
	\hspace*{0.02in}{\bf Step 3:} Get the unique solution $Q^{g_n}_{G_*}=(q_i^{g_n}:\ i \in G_*)$ of
		\begin{eqnarray}\label{5.3}				Q^{g_n}_{G_*}=(I-\emph{\textbf{P}}^{g_n}_{G_*})^{-1}V_B^{g_n}.
		\end{eqnarray}
		
	\hspace*{0.02in}{\bf Step 4:} Choose $g_{n+1}\in \Pi_d^s(F_*)$ such that
		\begin{eqnarray*}
			\begin{cases}
				g_{n+1}(i)\in \mathop{\rm{argmin}}\limits_{a\in A_{g_n}(i)}\{p(B|i,a)\\
				\quad \quad \quad \quad \quad+\sum\limits_{j \in G_*}p(j|i, a)
				q_j^{g_n}\},&\text{if}\ i\in G_*, A_{g_n}(i)\neq \emptyset\\
				g_{n+1}(i)=g_n(i), & \text{otherwise}.
			\end{cases}
		\end{eqnarray*}
		
	\hspace*{0.02in}{\bf Step 5:} A. If $g_{n+1}=g_n$, then stop, and set the optimal policy $g_*=g_n$ and the minimal reaching probability
		\begin{eqnarray*}
			q_i^*=\begin{cases}
				q_i^{g_n},& \text{if}\  i \in G_*\\
				0,\quad & \text{if}\  i \in F_*;
			\end{cases}
		\end{eqnarray*}
	
	B.	otherwise increment $n$ by $1$ and return to step 3.

\quad

The following proposition guarantees that  Algorithm 2 terminates after a finite number of iterations.

\begin{prop}\label{th5.3}
	For the state-classification-based PI algorithm
 in Algorithm 2,
	there exists $N<\infty$ such that
	\begin{eqnarray*}
		\nonumber	&&p(B|i,g_{_N}(i))+\sum_{j \in G_*}p(j|i, g_{_N}(i))q_j^{g_{_N}}\\
		&=&\min_{a\in A^*(i)}\{p(B|i,a)+\sum_{j \in G_*}p(j|i, a)q_j^{g_{_N}}\},\ i\in G_*.
	\end{eqnarray*}
	Hence, this policy $g_{_N}$ is optimal (i.e., $q_i^{g_{_N}}=q^*_i$ for $i\in B^c$).
\end{prop}

\begin{IEEEproof}
Take $g_n, g_{n+1} \in \Pi_d^s(F_*)$, which are the successive policies from  Algorithm 2. Then, $q_i^{g_{n+1}} \leq q_i^{g_n}$ for all $i \in B^c$ and $q^{g_{n+1}}_i<q_i^{g_n}$ for some $i \in B^c$. Since   $\Pi_d^s(F_*) \subset \Pi_d^s$ is finite, the algorithm must terminate under the stopping criterion in step 5, in a finite number $N$ of iterations. Thus, we have
\begin{eqnarray*}
		q_i^{g_{_N}}&=&p(B|i, g_{_{N}}(i))+\sum_{j \in G_*}p(j|i, g_{_{N}}(i))q_j^{g_{_N}}\\
		&=&\min_{a \in A_*(i)}\{ p(B|i, a)+\sum_{j \in G_*}p(j|i, a)q_j^{g_{_N}} \},\ \ i \in G_*,
\end{eqnarray*}
which, together with Theorem \ref{le5.1}, completes the proof.
\end{IEEEproof}
 Combining with Algorithm 1 and Algorithm 2, we get a new method of computing minimal reaching probability and its optimal policies, called the state-classification-based PI approach.
	
Since the complexity of different PI algorithms is still open (see \cite{L-95}), it would be important to analyze the number of calculations in each iteration of our state-classification-based PI approach, and compare it with that in the  PI algorithm for the  equivalent average  MDPs in \cite{DM22}. Here the calculations consist of addition, subtraction, multiplication and division. For this end, let
\begin{eqnarray*}
	 |G_*|=:\Gamma_1, \quad |B^c|=:\Gamma_2 \quad \text{and} \quad \max\limits_{i \in S}|A(i)|=:\Gamma_3.
\end{eqnarray*}

\begin{thm}\label{re-5.1}
\rm{ (On the effectiveness of Algorithms 1-2.) 
\begin{itemize}
\item[(i)]  The maximal total number of calculations for Algorithm 1 is given by
$N_1:=\frac{\Gamma_2\Gamma_3(\Gamma_2+1)}{2}$.
\item[(ii)]  The maximal number of calculations in one iteration in Algorithms~2, is
	$N_2=\frac{2\Gamma_1^3+15\Gamma_1^2-11\Gamma_1}{6}+2\Gamma_3\Gamma_1^2$.
\item[(iii)]  For the PI algorithm for the equivalent average MDPs in  \cite{DM22}, the maximal number of calculations in one iteration is 
$	N_3:=\frac{4\Gamma_2^3+30\Gamma_2^2-22\Gamma_2}{6}+4\Gamma_3\Gamma_2^2$.
\item[(iv)] 
$
N_3-(N_1+N_2) \geq \frac{2\Gamma_2^3+15\Gamma_2^2-11\Gamma_2}{6}+4\Gamma_3(\Gamma_2-\Gamma_1)
    >0.
$
\end{itemize}
}
\end{thm}	
	
\begin{IEEEproof}	
(i)\ For Algorithm 1, to obtain $U_n$ from $U_{n-1}$ at $n$'th iteration, we need to check whether, $p(\cup_{k=0}^{n-1} U_k|i, a)>0$ for every $i\in U_0^c \! \setminus\!\cup_{k=1}^{n-1}U_k$ and $a \in A(i)$ or not. At the $n$'th iteration, since $|U_k| \leq 1$ for all $k<n$, we can find out at least one state that belongs to $U_n$. Then, at the $n$'th iteration, because there are $\Gamma_3$ actions, the number of calculations is at most $(\Gamma_2-n+1)\Gamma_3$. Therefore, the total number of calculations is $N_1=\sum_{n=1}^{\Gamma_2-1}(\Gamma_2-n+1)\Gamma_3=\frac{\Gamma_3\Gamma_2(\Gamma_2+1)}{2}$.

(ii)\ (1)\ At step 3 in Algorithm 2, we use Gaussian elimination method to solve (\ref{5.2}) with $\Gamma_1$ variables, it needs to calculate $\frac{2\Gamma_1^3+15\Gamma_1^2-11\Gamma_1}{6}$ times.

(2)\ For step 4 in Algorithm 2, we need to calculate $p(B|i,a)+\sum_{j \in G_*}p(j|i, a)q_j^{g_n}$ with  $\Gamma_1$ equations  under a fixed $a$ and thus such number of calculations is
$2\Gamma_1$. Therefore, to obtain $g_{n+1}(i)$, due to $|A_{g_n}(i)|\leq \Gamma_3$ for every $i \in G_*$ and $g_n \in \Pi_d^s(F_*)$, we need to compute $p(B|i,a)+\sum_{j \in G_*}p(j|i, a)q_j^{g_n}$ for $a \in A_{g_n}(i)$ at most $\Gamma_3$ times. Therefore, for step 4, we have to calculate $2\Gamma_3\Gamma_1^2$ times.
Hence, for Algorithm 2, the number of calculations in one iteration is $
N_2=\frac{2\Gamma_1^3+15\Gamma_1^2-11\Gamma_1}{6}+2\Gamma_3\Gamma_1^2$.

(iii)\ Similarly, for the PI algorithm in \cite{DM22} and \cite{MLT05}, the number of calculations in one iteration is given as below.

(1)\ For the step of value approximation, we also use Gaussian elimination method to solve two systems of equations (9.2.1)-(9.2.2) in P. 452 from \cite{MLT05} with $\Gamma_2$ variables for each system of equations. Then, it needs to calculate $2\times\frac{2\Gamma_2^3+15\Gamma_2^2-11\Gamma_2}{6}$ times
in one iteration.

(2)\ For two steps of policy improvement under two different phases with $\Gamma_2$ variables, we need to calculate $2\Gamma_3\Gamma_2^2$ times in one iteration for each policy improvement.
Therefore, for the multi-chain policy iteration algorithm in \cite{MLT05}, the number of calculations in one iteration is given by $
N_3=\frac{4\Gamma_2^3+30\Gamma_2^2-22\Gamma_2}{6}+4\Gamma_3\Gamma_2^2$.

(iv)\ Finally, we have
\begin{eqnarray*}\label{L}
N_3-(N_1+N_2)\geq \frac{2\Gamma_2^3+15\Gamma_2^2-11\Gamma_2}{6}+4\Gamma_3(\Gamma_2-\Gamma_1)>0,
\end{eqnarray*}
where the first inequality is due to the fact that $\Gamma_2 > \Gamma_1$.
\end{IEEEproof}
\begin{rem}\label{re-eff}
	\rm{
		From Theorem \ref{re-5.1}(iv), according to our state-classification-based PI approach, the sum of total number of calculations in Algorithm 1 and the number of calculations in one iteration of Algorithm 2, is less than the number of calculations in one iteration of the PI algorithm for the equivalent average MDPs in \cite{DM22}. The larger the $\Gamma_2$ is, the obvious the effectiveness of our state-classification-based PI approach is.

}
\end{rem}

\section{Application to maintenance problems}
In \cite{A07,L96,DH08,N05,MYK04}, there are a large number of maintenance problems about the reliability engineering. This section applies our state-classification-based PI approach to a maintenance problem, in which the minimal reaching probability and an optimal policy are obtained by using Algorithm 2. Below we call the reaching probability as the failure probability.

\begin{exm}
	\rm{	(A maintenance problem)
		Consider a system consisting of two machines $M_1,\ M_2$ and a repairman. Each machine has three situations: available (labeled as $``2"$), deteriorates (labeled as $``1"$), and breaks down (labeled as $``0"$). An available machine, which does not need to be repaired, deteriorates with probability $\alpha_0>0$, breaks down with probability $\alpha_1>0$, and is still available with probability $1-\alpha_0-\alpha_1>0$.   We use a two-dimensional vector, $(i_1,i_2)$ $(i_1,i_2=0,1,2)$, to represent a state of the system, in which $i_k(k=1,2)$ denotes situations of machine $M_k$. For the analysis of our  optimization problem, we consider the following case: 1) When a machine breaks down, the machine can not be repaired by the repairman. 2) When a machine deteriorates, the repairman will repair it. 3) If there is at least one machine breaks down or both two machines deteriorates, the system must stop for checking. Thus, the system is working if and only if that both two machines are available or that one machine is available and another machine deteriorates.
		
		When one machine is available and another machine deteriorates, suppose that the repairman has two classes of ways (denoted by $c$ and $d$)  to deal with the deteriorated machine. If the repairman takes $c$, then the machine is available with probability $\beta_0$, still deteriorates with probability $\beta_1$, and breaks down with probability $1-\beta_0-\beta_1$. While if the repairman takes $d$, then the machine is available with probability $\theta_0$, still deteriorates with probability $\theta_1$, and breaks down with probability $1-\theta_0-\theta_1$. When the system stops for checking or has two machines which are available, the repairman does nothing.
		
		From the description of above problem, the state space of the system is $S=\{(i_1,i_2):\ i_1,i_2=0,1,2\}$. For convenience, set
		\begin{eqnarray*}
			S_1:=\{(0,0),(0,1),(0,2),(1,0),(1,1),(2,0)\}
		\end{eqnarray*}
		and when the system state is in $S_1$, the system stops for checking.  Then, $A(i_1,i_2)=\{c(i_1,i_2),d(i_1,i_2)\}$ for $(i_1,i_2) \in \{(1,2),(2,1)\}$, and $A(i_1,i_2)=\{\Delta\}$ for $(i_1,i_2)\in S\setminus \{(1,2),(2,1)\}$, where $\Delta$ means the action that the repairman does nothing.
		
		Obviously, the state $(0,0)$ is the worst one. Our goal is to find a policy under which the failure probability of the system is minimal. For this purpose, denote $B:=\{(0,0)\}$. Then, we only need to find an optimal policy $g_*$ such that the failure probability probability $P_{(i_1,i_2)}^{g_*}(\tau_{B}<\infty)$ is minimal.
		On the other hand, from the description of the maintenance problem, we see that the transition probabilities of the system are as follows:
		
		For $(i_1,i_2)\in S_1$,
		\begin{eqnarray*}
			\begin{cases}&p((i_1,i_2)|(i_1,i_2), \Delta)=1\\
				&p((j_1,j_2)|(i_1,i_2), \Delta)=0,\ \text{if}\ (j_1,j_2)\neq (i_1,i_2).
			\end{cases}
		\end{eqnarray*}
	
	For $(i_1,i_2)=(1,2)$,
		\begin{eqnarray*}\!
			\begin{cases}&p((0,0)|(1,2),a)\!=\!(1-\!\gamma_0\!-\!\gamma_1)\alpha_1\\
				&p((0,1)|(1,2),a)\!=\!(1-\!\gamma_0-\!\gamma_1)\alpha_0\\			&p((0,2)|(1,2),a)\!=\!(1\!-\!\gamma_0-\!\gamma_1)
				(1\!-\!\alpha_0-\!\alpha_1)\\
				&p((1,0)|(1,2),a)\!=\!\gamma_1 \alpha_1\\
				&p((1,1)|(1,2),a)\!=\!\gamma_1 \alpha_0\\
				&p((1,2)|(1,2),a)\!=\!\gamma_1(1\!\!-\!\alpha_0\!\!-\!\alpha_1)\\
				&p((2,0)|(1,2),a)\!=\!\gamma_0 \alpha_1\\
				&p((2,1)|(1,2),a)\!=\!\gamma_0 \alpha_0\\
				&p((2,2)|(1,2),a)\!=\!\gamma_0(1\!\!-\!\alpha_0\!\!-\!\alpha_1),
			\end{cases}
		\end{eqnarray*}
	where $(\gamma_0, \gamma_1)=(\beta_0,\beta_1)$ for  $a=c(1,2)$, and $(\gamma_0, \gamma_1)=(\theta_0,\theta_1)$ for $a=d(1,2)$.

	For $(i_1,i_2)=(2,1)$,
	\begin{eqnarray*}
		\begin{cases}&p((0,0)|(2,1),a)\!=\!\alpha_1(1\!-\!\gamma_0-\!\gamma_1)\\
			&p((0,1)|(2,1),a)\!=\!\alpha_1\gamma_1\\
			&p((0,2)|(2,1),a)\!=\!\alpha_1\gamma_0\\
			&p((1,0)|(2,1),a)\!=\!\alpha_0(1\!-\!\gamma_0-\!\gamma_1)\\
			&p((1,1)|(2,1),a)\!=\!\alpha_0\gamma_1\\
			&p((1,2)|(2,1),a)\!=\!\alpha_0\gamma_0\\			&p((2,0)|(2,1),a)\!=\!(1\!-\!\alpha_0-\!\alpha_1)
			(1\!-\!\gamma_0-\!\gamma_1)\\
			&p((2,1)|(2,1),a)\!=\!(1-\!\alpha_0-\!\alpha_1)\gamma_1\\
			&p((2,2)|(2,1),a)\!=\!(1-\!\alpha_0-\!\alpha_1)\gamma_0,
		\end{cases}
	\end{eqnarray*}
	where $(\gamma_0, \gamma_1)=(\beta_0,\beta_1)$ for $a=c(2,1)$, and $(\gamma_0, \gamma_1)=(\theta_0,\theta_1)$ for $a=d(2,1)$.

	For $(i_1,i_2)=(2,2)$,
	\begin{eqnarray*}
		&&p((j_1,j_2)|(2,2),\Delta)\\
		&=&\begin{cases}
			\alpha_1^2,\quad \quad \quad \quad \quad \quad & \ \text{if} \  (j_1,j_2)=(0,0)\\
			\alpha_1\alpha_0,\quad \quad \quad \quad \quad & \ \text{if} \  (j_1,j_2)=(0,1)\\
			\alpha_1(1-\alpha_0-\alpha_1),\ & \ \text{if} \  (j_1,j_2)=(0,2)\\
			\alpha_0\alpha_1,\quad \ \ \  \quad \quad & \ \text{if} \  (j_1,j_2)=(1,0)\\
			\alpha_0^2,\ \quad \quad \quad \quad \quad & \ \text{if} \  (j_1,j_2)=(1,1)\\
			\alpha_0(1-\alpha_0-\alpha_1),\ & \ \text{if} \  (j_1,j_2)=(1,2)\\
			(1-\alpha_0-\alpha_1)\alpha_1,\ & \ \text{if} \  (j_1,j_2)=(2,0)\\
			(1-\alpha_0-\alpha_1)\alpha_0,\ & \ \text{if} \  (j_1,j_2)=(2,1)\\
			(1-\alpha_0-\alpha_1)^2,\ & \ \text{if}  \  (j_1,j_2)=(2,2).
		\end{cases}
	\end{eqnarray*}
We now use our state-classification-based PI approach to compute the minimal failure probability as below. By Algorithm~1, we have
	\begin{eqnarray*}
		F_*=\{(0,1),(0,2),(1,0),(1,1),(2,0)\},\ G_*=\{(1,2),(2,1),(2,2)\}
	\end{eqnarray*}
	and $\Pi_d^s(F_*)=\{f_1, f_2, f_3, f_4\}$, where
	\begin{eqnarray*}
	f_1(i_1,i_2)=\begin{cases}
		c(i_1,i_2),\ & \text{if}\ (i_1,i_2) \in \{(1,2),(2,1)\},\\
		\Delta,\ &  \text{otherwise},
	\end{cases}
\end{eqnarray*}
\begin{eqnarray*}
	f_2(i_1,i_2)=\begin{cases}
		c(1,2),\ & \text{if}\ (i_1,i_2)=(1,2),\\
		d(2,1),\ & \text{if}\ (i_1,i_2)=(2,1),\\
		\Delta,\ &  \text{otherwise},
	\end{cases}
\end{eqnarray*}
\begin{eqnarray*}
	f_3(i_1,i_2)=		\begin{cases}d(i_1,i_2),\ & \text{if}\ (i_1,i_2) \in \{(1,2),(2,1)\},\\
		\Delta,\ &  \text{otherwise},
	\end{cases}
\end{eqnarray*}
\begin{eqnarray*}
	f_4(i_1,i_2)=\begin{cases}
		d(1,2),\ & \text{if}\ (i_1,i_2)=(1,2),\\
		c(2,1),\ & \text{if}\ (i_1,i_2)=(2,1),\\
		\Delta,\ & \text{otherwise}.
	\end{cases}
\end{eqnarray*}

	For numerical computation, we take $\alpha_0=\alpha_1=\frac{1}{4}$.
	
	Applying Algorithm 2, we finally get that if $``5\beta_0+6\beta_1\geq 5\theta_0+6\theta_1"$, then $f_1$ is optimal and the minimal failure probability is
	\begin{eqnarray*}
		\begin{cases}
			&q^*_{(1,2)}=\frac{6-5\beta_0-6\beta_1}{2(12-5\beta_0-6\beta_1)}\\
			&q^*_{(2,1)}=\frac{6-5\beta_0-6\beta_1}{2(12-5\beta_0-6\beta_1)}\\
			&q^*_{(2,2)}=\frac{24-15\beta_0-18\beta_1}{12(12-5\beta_0-6\beta_1)}.
		\end{cases}
	\end{eqnarray*}
	Similarly, for the other case of $``5\beta_0+6\beta_1<5\theta_0+6\theta_1"$, we know that $f_3$ is optimal and the minimal failure probability is
	\begin{eqnarray*}
		\begin{cases}
			&q^*_{(1,2)}=\frac{6-5\theta_0-6\theta_1}{2(12-5\theta_0-6\theta_1)}\\		&q^*_{(2,1)}=\frac{6-5\theta_0-6\theta_1}{2(12-5\theta_0-6\theta_1)}\\		 &q^*_{(2,2)}=\frac{24-15\theta_0-18\theta_1}{12(12-5\theta_0-6\theta_1)}.
		\end{cases}
	\end{eqnarray*}

}

\end{exm}

\section{Conclusion}
For MDPs with finite state and action spaces, a new method for computing the minimal reaching probability and its optimal policy, is our paper's main contribution, and the detail is illustrated as below.
\begin{itemize}
	\item[(1)] After introducing  the concept of an absorbing set of a stationary policy, we have presented a $B^c$-state-classification method to calculate the largest absorbing set in $B^c$ (i.e., the Algorithm 1), which can be used for  establishing a new IOE with a unique solution, that is, the minimal reaching probability.
	\item[(2)] By the uniqueness of solution of the IOE and the $B^c$-state-classification method, a state-classification-based PI algorithm (i.e., the Algorithm 2) of optimal policies and the minimal reaching probability, has been obtained.
	
	\item[(3)] Combining with (1) and (2), our state-classification-based PI approach is obtained. By comparing it with the existing  PI algorithm for the equivalent average MDPs in \cite{DM22}, we find that the sum of total number of calculations in Algorithm 1 and the number of calculations in one iteration of Algorithm 2, is less than the number of calculations in one iteration of the PI algorithm for  the equivalent average MDPs in \cite{DM22}.
\end{itemize}

\section*{Declaration of competing interest}
The authors declare that they have no known competing financial interests or personal relationships that could have appeared to influence the work reported in this paper.

\section*{Data availability}
No data was used for the research described in the article.

\section*{Acknowledgment}


\ifCLASSOPTIONcaptionsoff
\newpage
\fi
This work is supported by the National Key Research and Development Program of China (2022YFA1004600).

\end{document}